\documentclass[11pt]{elsarticle}
\usepackage{amsmath,amssymb,amsthm}
\usepackage{hyperref}
\usepackage{mathrsfs}
\usepackage{graphicx}
\usepackage{tikz}
\usepackage{enumerate}
\usepackage{float}

\theoremstyle{plain}
\newtheorem{theorem}{Theorem}[section]
\newtheorem{lemma}[theorem]{Lemma}

\newtheorem{proposition}[theorem]{Proposition}

\theoremstyle{definition}

\newtheorem{remark}[theorem]{Remark}

\newtheorem{problem}[theorem]{Problem}

\newcommand{\ua}{\mathord{\uparrow}}
\newcommand{\da}{\mathord{\downarrow}}
\newcommand{\rom}[1]{\rm{\uppercase\expandafter{\romannumeral #1}}}

\makeatletter
\def\ps@pprintTitle{%
  \let\@oddhead\@empty
  \let\@evenhead\@empty
  \def\@oddfoot{\reset@font\hfil\thepage\hfil}
  \let\@evenfoot\@oddfoot
}
\makeatother

\begin{document}

\begin{frontmatter}

\title{The core compactly generated topology \tnoteref{t1}}
\tnotetext[t1]{This work is supported by the National Natural Science Foundation of China (No.12231007) and by Hunan Provincial Innovation Foundation For Postgraduate (CX20200419)}

\author{Qingguo Li}
\ead{liqingguoli@aliyun.com}
\author{Hualin Miao\corref{a1}}
\ead{miaohualinmiao@163.com}
\address{School of Mathematics, Hunan University, Changsha, Hunan, 410082, China}
\cortext[a1]{Corresponding author.}
\begin{abstract}
M. Escard$\mathrm{\acute{o}}$ et al. asked whether the core compactly generated topology of a sober space is again sober and the sobrification of a core compactly generated space again core compactly generated. In this note, we answer the problem by displaying a counterexample, which reveals that the core compactly generated spaces are not closed under sobrifications. Meantime, we obtain that the core compactly generated spaces are closed under $\omega$-well-filterifications and $D$-completions. Furthermore, we find that the core compactly generated topology of the Smyth power space of a well-filtered space coincides with the Scott topology. Finally, we provide a characterization for the core-compactness of core compactly generated spaces.
\end{abstract}

\end{frontmatter}
\section{Introduction}
 Sobriety, well-filteredness and monotone convergence have been motivating a widespread concern and being the most important and useful properties in Non-Hausdorff topology. Moreover, they have been extensively investigated in \cite{Ershov},\cite{GHKLMS},\cite{GL},\cite{Keimel},\cite{Miao1},\cite{Miao2},\cite{Miao3},\cite{SXXZ},\cite{Wyler},\cite{WXXZ},\cite{XSXZ},\cite{zhongxizhang},\cite{ZF} and so on.

In \cite{Day}, Day considered an enlarged category of topological spaces which is cartesian closed, and showed this category can be reflected onto the category of $\mathcal{C}$-generated spaces. Based on this fact and the general categorical reflection theorem, he deduced that the category of $\mathcal{C}$ generated spaces is cartesian closed.

Afterwards, M. Escard$\mathrm{\acute{o}}$ et al. presented a simple and uniform proof of Cartesian closedness for the category of $\mathcal{C}$-generated spaces, including the category of compactly generated spaces and core compactly generated spaces (\cite{Escardo}).

In \cite{Battenbeld}, O. Battenfeld et al. presented the result that the $D$-completion of a compactly generated space is again compactly generated. Furthermore, G. Gruenhage and T. Streicher reached the conclusion that the sobrification of a compactly generated space may not be compactly generated \cite{Gruenhage}. Other scholars also considered the category of compactly generated spaces (\cite{M. Escardo}, \cite{Vogt}).

As for the core compactly generated spaces, M. Escard$\mathrm{\acute{o}}$ et al. posed the following problem \cite{Escardo}:
\begin{problem}
Is the core compactly generated topology of a sober space again sober?
Dually, is the sobrification of a core compactly generated space again core compactly
generated?
\end{problem}
 In this note, we give a counterexample to illustrate that the core compactly generated topology of a sober space may not be sober and the sobrification of a core compactly generated space may not be core compactly generated. However, we discover that the core compactly generated topology of a $\omega$-well-filtered space (resp. a monotone convergence space) is again $\omega$-well-filtered (resp. monotone convergence) and the $\omega$-well-filterification (resp. $D$-completion) of a core compactly generated space is also core compactly generated. Moreover, we find that the core compactly generated topology of a monotone convergence space is contained in the Scott topology. Naturally, we wonder on which class of topological spaces, the core compactly generated topologies are exactly the Scott topologies determined by the specialization orders on those spaces. In Section 5, we confirm that the Smyth power space of each well-filtered space can satisfy that its core compactly generated topology and Scott topology coincide.

 In \cite{Escardo}, M. Escard$\mathrm{\acute{o}}$ et al. derived that if a topology space $X$ is core compact, then the topological product of $X,Y$ is also core compactly generated for any core compactly generated space $Y$. We conclude the paper by showing that a core compactly generated space $X$ is core compact iff the topological product of $X,Y$ is also core compactly generated for any core compactly generated space $Y$.
\section{Preliminaries}
Let $P$ be a poset. A subset $D$ of $P$ is \emph{directed} provided it is nonempty and every finite subset of $D$ has an upper bound in $D$. A poset $P$ is a \emph{dcpo} if every directed subset $D$ has a supremum. A subset $U$ of a poset $P$ is \emph{Scott open} if (i) it is an upper set $(U = \ua U = \{x \in P : u \leq x ~\mathrm{for ~some} ~u \in U\})$ and (ii) for every
directed subset $D$ of $P$ with $\sup D$ existing and $\sup D\in U$, it follows that $D \cap U \neq \emptyset$. The complement of a Scott open set is called \emph{Scott closed}. Let $\sigma (P)$ denote the set of all Scott open sets and $\Gamma(P)$ the set of all Scott closed sets. The space $(P,\sigma(P))$ called the Scott space of $P$ is written as $\Sigma P$. Without further references, the posets mentioned here are all endowed with the Scott topology.

 Given a $T_{0}$ space $(X, \tau)$, we denote the closure of $A$ by $cl_{\tau}(A)$ for any subset $ A\subseteq X$. We define $x\leq y$ iff $x\in cl(y)$. Hence, $X$ with its specialization order is a poset. We denote the set of all open sets by $\mathcal{O} (X)$. The set $A$ is called \emph{$d$-closed} if $D\subseteq A$ implies that $\sup D\in A$ for any directed subset $D$ of $A$. Let $cl_{d}(A)$ denote the closure of $A$ in the $d$-topology. The set $A$ is called a \emph{$d$-dense} subset of $X$ if the $d$-closure of $A$ is $X$, i.e., $cl_{d} (A)=X$. The set of all irreducible closed subsets of $X$ is denoted by $\mathbf{IRR}(X)$. A set $K$ of a topological space is called \emph{saturated} if it is the intersection of all open sets containing $K$ ($K=\ua K$ in its specialization order). For a topological space $X$, the set of all non-empty compact saturated subsets of $X$ are denoted by $Q(X)$. 

\section{The core compactly generated topology of sober spaces}

 From \cite{Escardo}, we have that the core compactly generated topology of a $T_{0}$ space $X$ enjoys the same specialization order as $X$ and a dcpo endowed with the Scott topology is core compactly generated.

 The goal of this section is to give a negative answer to the problem posed by M. Escard$\mathrm{\acute{o}}$ et al.~\cite{Escardo} mentioned in the introduction.
 \begin{remark}\label{transfinite}
   Let $L$ be a dcpo. If an upper set $U$ of $L$ has the property, that $\sup C\in U$ suggests that $C\cap U\neq \emptyset$, for any well-ordered chain $C$ of $L$, then $U$ is Scott open.
   \begin{proof}
    Let $U$ be an upper set, with the property that $\sup C\in U$ indicates that $C\cap U\neq \emptyset$ for any well-ordered chain $C$ of $L$. We use transfinite induction on the cardinality of the directed
set $D$.

If $D$ is finite, and $\sup D\in U$, then $\sup D\in D\cap U$.

Now assume that $D$ is infinite and $E\cap U\neq \emptyset$ for any directed set $E$ with cardinality smaller than ${\mid }D{\mid}$ and $\sup E\in U$. From the theorem of Iwomura \cite{Markowsky}, $\sup D=\sup_{\alpha<{\mid}D{\mid}}\sup D_{\alpha}$ and $\{\sup D_{\alpha}\mid \alpha<{\mid}D{\mid}\}$ is a well-ordered chain with $D_{\alpha}\subseteq D$, which yield that there is $\alpha< {\mid}D{\mid}$ such that $\sup D_{\alpha}\in U$. By induction hypothesis, $D_{\alpha}\cap U\neq \emptyset$, whence $D\cap U\neq \emptyset$.
   \end{proof}
 \end{remark}
 Recall that a $T_0$ space is a \emph{monotone convergence space} if and only if the closure of every directed set (in the specialization order) is the closure of a unique point.
 A subset $V$ of a topological space $(X,\tau)$ is \emph{open in the core compactly generated topology} on $X$ if, for every core compact space $C$ and continuous function $\rho: C \rightarrow X$, the preimage $p ^{-1} (V)$ is open in $C$. We write $\mathcal{C}X$ for
the core compactly generated topology of $X$, and we say that $X$ is \emph{core compactly generated} if $\tau =\mathcal{C}X$. Note that it is direct that $\tau \subseteq \mathcal{C}X$ by the above concept.
 \begin{lemma}\label{Scott}
   Let $(X,\tau)$ be a monotone convergence space. Then the core compactly generated topology $\mathcal{C}X$ of $X$ is contained in the Scott topology of $X$ under the specialization order of $X$.
   \begin{proof}
     Due to \cite[Lemma~4.6]{Escardo}, we know that $U$ is an upper set for any $U\in \mathcal{C}X$. Now let $C$ be a well-ordered chain of $X$ with $\sup C\in U$. Define $C^{'}=\{\sup E\mid E\subseteq C\}$. One sees immediately that $C^{'}$ is a subdcpo of $X$ and $C^{'}$ is a complete chain. This reveals that the inclusion map $i: \Sigma C^{'}\rightarrow (X,\tau)$ is continuous in the light of the monotone convergence of $X$. Note that $\Sigma C^{'}$ is core compact. We conclude that $i^{-1}(U)$ is Scott open in $C^{'}$ since $U$ is core compactly generated. The fact that $\sup C\in i^{-1}(U)$ deduces that $C\cap U\neq \emptyset$ with the help of $C\subseteq C^{'}$. By applying Remark \ref{transfinite}, the set $U$ is Scott open.
   \end{proof}
 \end{lemma}
 We shall say that a $T_{0}$ space $X$ is \emph{well-filtered} if for each filter basis $\mathcal{C}$ of compact saturated sets and each open set $U$ with $\bigcap \mathcal{C} \subseteq U$, there is a $K\in \mathcal{C}$ with $K\subseteq U$.
  An arbitrary nonempty subset $A$ of a $T_{0}$ space $X$ is \emph{irreducible} if $A \subseteq  B\cup C$ for closed subsets $B$ and $C$ implies $A\subseteq B$ or $A \subseteq C$. A topological space $X$ is \emph{sober} if it is $T_{0}$ and every irreducible closed subset of $X$ is the closure of a (unique) point.
A sobrification of a $T_{0}$ space $X$ consists of a sober space $Y$ and a continuous map $\eta : X \rightarrow Y$ which enjoys the following universal property: For every continuous map $f$ from $X$ to a sober space $Z$, there is a unique continuous map $\overline{f}: Y \rightarrow Z$ such that $f=\overline{f} \circ \eta$.

A standard construction for the sobrification of a $T_{0}$ space $X$ is to set

\begin{center}
  $X^{s} := \{ A\subseteq X : A\in \mathbf{IRR}(X)\}$
\end{center}
topologized by open sets $U^{ s}:=\{ A \in X ^{s} : A \cap U \neq \emptyset\}$ for each open subset $U$ of $X$. If we define $\eta^{s}_{X} :X \rightarrow X ^{s}$ by $\eta^{s}_{X} ( x ) = cl(\{ x \}) $, then we obtain a sobrification \cite[Exercise~V-4.9]{GHKLMS}, which we call the standard sobrification.

\begin{proposition}\label{core}
  Let $\Sigma P$ be a well-filtered dcpo which is not sober and $\mathbf{IRR}(P)=\{\da x\mid x\in P\}\cup\{P\}$. Then the core compactly generated topology of the standard sobrification $P^{s}$ of $\Sigma P$ is the Scott topology of $(\mathbf{IRR}(P),\subseteq)$.
  \begin{proof}
   Through Lemma \ref{Scott}, it remains to prove that the Scott topology of $\mathbf{IRR}(P)$ is contained in the core compactly generated topology. Obviously, we see that $\sigma(\mathbf{IRR}(P))=\{U^{s}\mid U\in \sigma(P)\}\cup \{\{P\}\}$. So it suffices to show that $\{P\}$ is core compactly generated open. By \cite[Lemma 8.2]{Escardo}, we just need to check that $\{P\}$ is locally compact sober generated open.

   To this end, let $\rho: C\rightarrow P^{s}$ be a continuous function from a locally compact sober space $C$. For any $x\in \rho^{-1}(\{P\})$. The fact that $C$ is locally compact implies that $\ua x=\bigcap _{x\in int(K),K\in Q(C)}K$. On account of \cite[Lemma 8.1]{Escardo}, we derive that $\bigcap_{x\in int(K),K\in Q(C)}\ua \rho(K)=\ua \rho(\ua x)=\{P\}$. Due to the continuity of $\rho$, we demonstrate that $\{\ua \rho(K)\mid x\in int(K),K\in Q(C)\}$ is a filtered family of $Q(P^{s})$.

    The assumption that $\{P\}$ is the greatest element of $\mathbf{IRR}(P)$ guarantees that $\{\ua \rho(K)\backslash\{P\}\mid x\in int(K),K\in Q(C)\}$ is a filtered family of $Q(P^{s}\backslash\{P\})$ and $\bigcap_{x\in int(K),K\in Q(C)}(\ua \rho(K)\backslash\{P\})=\emptyset$. What is noteworthy is that $\Sigma P$ is homeomorphic to the subspace $P^{s}\backslash \{P\}$ of $P^{s}$. This means that $P^{s}\backslash\{P\}$ is well-filtered owing to the assumption that $\Sigma P$ is well-filtered. Hence, there exists $K\in Q(C)$ with $x\in int(K)$ such that $\ua \rho(K)\backslash\{P\}=\emptyset$. In other words, $x\in int(K)\subseteq K\subseteq \rho^{-1}(\{P\})$. Now we gain our desired result.
  \end{proof}
\end{proposition}

Let us consider the well-filtered dcpo $P=\mathbb{N}\times\mathbb{N}\times(\mathbb{N}\cup \{\infty\})$. The order  $(i_1,j_1,m_1)\leqslant (i_2,j_2,m_2)$ is defined as follows:

  $\bullet$ $i_1=i_2,j_1=j_2, m_1\leqslant m_2\leqslant \infty$;

  $\bullet$ $i_2=i_1+1,m_1\leqslant j_2$, $m_2=\infty$.

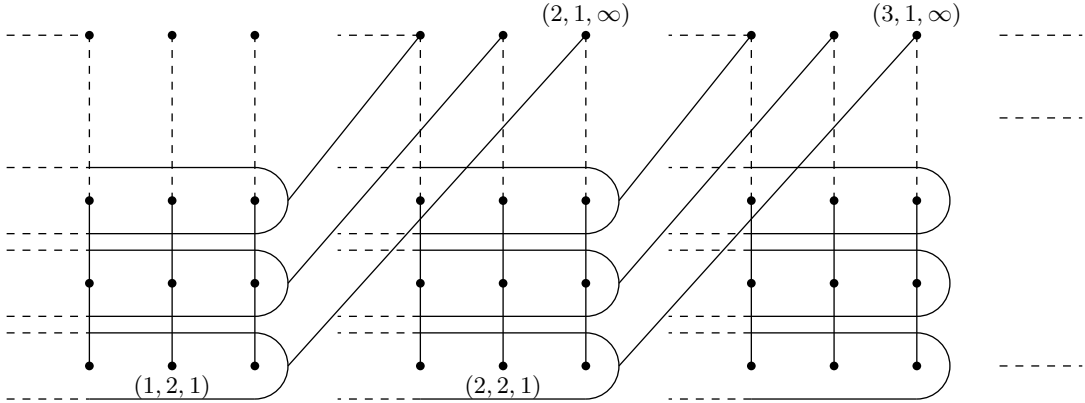
  \begin{tikzpicture}[line width=0.5pt,scale=1.1]
\fill[black] (0,0) circle (1.5pt);
\fill[black] (0,1) circle (1.5pt);
\fill[black] (0,2) circle (1.5pt);
\fill[black] (0,4) circle (1.5pt);
\draw (0,0)--(0,2);
\draw [dashed](0,2)--(0,4);
\fill[black] (1,0) circle (1.5pt);
\fill[black] (1,1) circle (1.5pt);
\fill[black] (1,2) circle (1.5pt);
\fill[black] (1,4) circle (1.5pt);
\draw (1,0)--(1,2);
\draw [dashed](1,2)--(1,4);
\fill[black] (2,0) circle (1.5pt);
\fill[black] (2,1) circle (1.5pt);
\fill[black] (2,2) circle (1.5pt);
\fill[black] (2,4) circle (1.5pt);
\draw (2,0)--(2,2);
\draw [dashed](2,2)--(2,4);
\draw [dashed](-1,4)--(0,4);
\draw (0,0.4)--(2,0.4);
\draw (0,-0.4)--(2,-0.4);
\draw (2,-0.4) arc (-90:90:0.4);
\draw [dashed](-1,0.4)--(0,0.4);
\draw [dashed](-1,-0.4)--(0,-0.4);

\draw (0,0.6)--(2,0.6);
\draw (0,1.4)--(2,1.4);
\draw (2,0.6) arc (-90:90:0.4);
\draw [dashed](-1,1.4)--(0,1.4);
\draw [dashed](-1,0.6)--(0,0.6);

\draw (0,1.6)--(2,1.6);
\draw (0,2.4)--(2,2.4);
\draw (2,1.6) arc (-90:90:0.4);
\draw [dashed](-1,2.4)--(0,2.4);
\draw [dashed](-1,1.6)--(0,1.6);


\fill[black] (6,0) circle (1.5pt);
\fill[black] (6,1) circle (1.5pt);
\fill[black] (6,2) circle (1.5pt);
\fill[black] (6,4) circle (1.5pt);
\draw (6,0)--(6,2);
\draw [dashed](6,2)--(6,4);
\fill[black] (4,0) circle (1.5pt);
\fill[black] (4,1) circle (1.5pt);
\fill[black] (4,2) circle (1.5pt);
\fill[black] (4,4) circle (1.5pt);
\draw (4,0)--(4,2);
\draw [dashed](4,2)--(4,4);
\fill[black] (5,0) circle (1.5pt);
\fill[black] (5,1) circle (1.5pt);
\fill[black] (5,2) circle (1.5pt);
\fill[black] (5,4) circle (1.5pt);
\draw (5,0)--(5,2);
\draw [dashed](5,2)--(5,4);
\draw [dashed](4,4)--(3,4);
\draw (4,0.4)--(6,0.4);
\draw (4,-0.4)--(6,-0.4);
\draw (6,-0.4) arc (-90:90:0.4);
\draw [dashed](4,0.4)--(3,0.4);
\draw [dashed](4,-0.4)--(3,-0.4);

\draw (4,0.6)--(6,0.6);
\draw (4,1.4)--(6,1.4);
\draw (6,0.6) arc (-90:90:0.4);
\draw [dashed](4,1.4)--(3,1.4);
\draw [dashed](4,0.6)--(3,0.6);

\draw (4,1.6)--(6,1.6);
\draw (4,2.4)--(6,2.4);
\draw (6,1.6) arc (-90:90:0.4);
\draw [dashed](4,2.4)--(3,2.4);
\draw [dashed](4,1.6)--(3,1.6);


\fill[black] (8,0) circle (1.5pt);
\fill[black] (8,1) circle (1.5pt);
\fill[black] (8,2) circle (1.5pt);
\fill[black] (8,4) circle (1.5pt);
\draw (8,0)--(8,2);
\draw [dashed](8,2)--(8,4);
\fill[black] (9,0) circle (1.5pt);
\fill[black] (9,1) circle (1.5pt);
\fill[black] (9,2) circle (1.5pt);
\fill[black] (9,4) circle (1.5pt);
\draw (9,0)--(9,2);
\draw [dashed](9,2)--(9,4);
\fill[black] (10,0) circle (1.5pt);
\fill[black] (10,1) circle (1.5pt);
\fill[black] (10,2) circle (1.5pt);
\fill[black] (10,4) circle (1.5pt);
\draw (10,0)--(10,2);
\draw [dashed](10,2)--(10,4);
\draw [dashed](8,4)--(7,4);
\draw (8,0.4)--(10,0.4);
\draw (8,-0.4)--(10,-0.4);
\draw (10,-0.4) arc (-90:90:0.4);
\draw [dashed](8,0.4)--(7,0.4);
\draw [dashed](8,-0.4)--(7,-0.4);

\draw (8,0.6)--(10,0.6);
\draw (8,1.4)--(10,1.4);
\draw (10,0.6) arc (-90:90:0.4);
\draw [dashed](8,1.4)--(7,1.4);
\draw [dashed](8,0.6)--(7,0.6);

\draw (8,1.6)--(10,1.6);
\draw (8,2.4)--(10,2.4);
\draw (10,1.6) arc (-90:90:0.4);
\draw [dashed](8,2.4)--(7,2.4);
\draw [dashed](8,1.6)--(7,1.6);

\draw (2.4,0)--(6,4);
\draw (2.4,1)--(5,4);
\draw (2.4,2)--(4,4);
\draw (6.4,0)--(10,4);
\draw (6.4,1)--(9,4);
\draw (6.4,2)--(8,4);

\draw [dashed](11,0)--(12,0);
\draw [dashed](11,4)--(12,4);
\draw [dashed](11,3)--(12,3);

\node[font=\footnotesize] (1) at(1,-0.25) {$(1,2,1)$};
\node[font=\footnotesize] (1) at(5,-0.25) {$(2,2,1)$};
\node[font=\footnotesize] (2) at(6,4.25) {$(2,1,\infty)$};
\node[font=\footnotesize] (3) at(10,4.25) {$(3,1,\infty)$};
\end{tikzpicture}
 \vskip 0.2cm

The above example is constructed by Jia \cite{Jia}, we can see that the dcpo $P$ is not sober.
\begin{theorem}
Let $P$ be the dcpo mentioned above. Then $P^{s}$ is a sober space. But the core compactly generated topology of $P^{s}$ is not sober. Moreover, $\Sigma P$ is core compactly generated, but $P^{s}$ is not core compactly generated.
\begin{proof}
  Due to Proposition \ref{core}, we can have that the core compactly generated topology of the standard sobrification $P^{s}$ of $\Sigma P$ is the Scott topology of $(\mathbf{IRR}(P),\subseteq)$. One sees clearly that $\Sigma P^{s}$ is not sober. Note that $\Sigma P$ is core compactly generated. Again  by applying Proposition \ref{core}, we know that $P^{s}$ is not core compactly generated.
\end{proof}
\end{theorem}
\section{The core compactly generated topology of $\omega$-well-filtered spaces}
From the above results, we conclude that the category of all core compactly generated sober spaces is not a reflective subcategory of the category of all core compactly generated spaces with continuous functions.

We now move on to the next goal of this note, which pays attention to the core compactly generated topology of $\omega$-well-filtered spaces. In this section, we prove that the core compactly generated topology of a $\omega$-well-filtered space is again $\omega$-well-filtered and the $\omega$-well-filterification of a core compactly generated space is again core compactly generated. This states that the category of all core compactly generated $\omega$-well-filtered spaces is a reflective subcategory of the category of all core compactly generated spaces with continuous functions.

  Recall that a $T_{0}$ space $X$ is called \emph{$\omega$-well-filtered}, if for any countable filtered family $\{K_{i} : i < \omega \}\subseteq Q(X)$ and $U \in \mathcal{O}(X)$, the following condition holds,
  \begin{center}
    $\bigcap_{i<\omega}K_{i}\subseteq U\Rightarrow \exists i_{0}< \omega, K_{i_{0}}\subseteq U$
  \end{center}

Let $X$ be a $T_{0}$ space. A nonempty subset $A$ of $X$ is said to
have the \emph{countable compactly filtered property} ($KF_{\omega}$ property), if there exists a filter family $\mathcal{K}$ of $Q(X)$ such that $cl(A)$ is a minimal closed set that intersects all members of $\mathcal{K}$, where $\mathcal{K}$ is a countable subset of $Q(X)$. We call such a set $KF_{\omega}$, or a $KF_{\omega}$-set. Denote by $\mathbf{KF_{\omega}}(X)$ the set of all closed
$KF_{\omega}$ subsets of $X$. The countable compactly filtered property is called $\omega$-Rudin property in \cite{XiaoquanXu}. 

\begin{lemma}\label{well-filtered}
  Let $(X,\tau)$ be a $\omega$-well-filtered space. Then the core compactly generated topology $\mathcal{C}X$ of $X$ is again $\omega$-well-filtered.
  \begin{proof}
    Assume for the sake of a contradiction that the core compactly generated topology $\mathcal{C}X$ of $X$ is not $\omega$-well-filtered. Via \cite[Corollary~6.11]{XiaoquanXu}, there is a closed $KF_{\omega}$-set $A$ in the core compactly generated topology with $A\notin\{\da x\mid x\in X\}$. It follows that there exists a countable descending chain $K_{0}\supseteq K_{1}\supseteq\cdots \supseteq K_{n}\supseteq \cdots $ of compact saturated subsets of $\mathcal{C}X$ such that $A$ is a minimal closed set in $\mathcal{C}X$ that intersects all members of $\{K_{n}\mid n\in \mathbb{N}\}$. Pick $x_{n}\in A\cap K_{n}$ for any $n\in \mathbb{N}$. We define $H=\{x_{n}\mid n\in \mathbb{N}\}$. Then we have that $cl_{\mathcal{C}X}(H)=A$ by the minimality of $A$.

     Consider the function $i:(X,\mathcal{C}X)\rightarrow (X,\tau)$ defined by $i(x)=x$ for any $x\in X$.
    Then the function $i$ is continuous, which implies that $A$ is a $KF_{\omega}$-set of $(X,\tau)$. Again applying Corollary 6.11 in \cite{XiaoquanXu}, we draw the conclusion that $\sup A\in cl_{\tau}(A)$ since $(X,\tau)$ is $\omega$-well-filtered.

    Claim 1: $\sup A\in cl_{\tau}(H)$.

    For any $U\in \tau$ with $\sup A\in U$, we can deduce that $U\cap A\neq \emptyset$ because of the fact that $\sup A\in cl_{\tau}(A)$. Choose $a\in U\cap A$. Note that $U\in \mathcal{C}X$ and $cl_{\mathcal{C}X}(H)=A$. Then $U\cap H\neq \emptyset$.

    Claim 2: If $U\cap H\neq \emptyset$, then $H\backslash U$ is finite for any $U\in \tau$.

    Now let $U\in \tau$ with $U\cap H\neq \emptyset$. Suppose $H\backslash U$ is infinite. Then $A\backslash U\cap K_{n}\neq \emptyset$ for any $n\in \mathbb{N}$. The fact that $U\in \tau \subseteq \mathcal{C}X$ suggests that $A\backslash U$ is closed in $\mathcal{C}X$. We conclude that $A=A\backslash U$ from the minimality of $A$. The assumption that $H\subseteq A$ indicates that $H\cap U=\emptyset$, which contradicts $H\cap U\neq \emptyset$. Denote the one-point compactification of the countable discrete space $\mathbb{N}$ by $\mathbb{N}^{\infty}$. Consider the function $\rho: N^{\infty}\rightarrow (X,\tau)$ defined by
     \begin{center}
  $$ \rho(x)=\left\{
\begin{aligned}
\sup A &,&x=\infty \\
x_{n} &, & x=n\in \mathbb{N} \\
\end{aligned}
\right.
$$
  \end{center}

  Claims 1 and 2 yield that $\rho$ is continuous. Then $\rho^{-1}(X\backslash A)$ is open in $N^{\infty}$ due to the generated openness of $X\backslash A$ and the core compactness of $\mathbb{N}^{\infty}$. One sees immediately that $N^{\infty}\in \rho^{-1}(X\backslash A)$. This implies that there is $n\in \mathbb{N}$ such that $n\in \rho^{-1}(X\backslash A) $, which is equivalent to saying that $x_{n}\in X\backslash A$. But $x_{n}\in A$, a contradiction.
  \end{proof}
\end{lemma}
Let us recall some known facts concerning the $\omega$-well-filterification of a $T_{0}$ space. A $\omega$-well-filterification of a $T_{0}$ space $X$ consists of a $\omega$-well-filtered space $Y$ and a continuous map $\eta : X \rightarrow Y$ which enjoys the following universal property: For every continuous map $f$ from $X$ to a $\omega$-well-filtered space $Z$, there is a unique continuous map $\overline{f}: Y \rightarrow Z$ such that $f=\overline{f} \circ \eta$.
A subset $A$ of a $T_{0}$ space $X$ is called a \emph{$\omega$-well-filtered determined set}, $WD_{\omega}$-set for short, if for any continuous mapping $f : X\rightarrow Y$ to a $\omega$-well-filtered space $Y$, there exists a unique $y_{A} \in Y$ such that $cl(f(A)) = cl(\{y_{A}\})$. The set of all closed $\omega$-well-filtered determined subsets of $X$ is denoted by $\mathbf{WD_{\omega}}(X)$.

A standard construction for the $\omega$-well-filterification of a $T_{0}$ space $X$ is to set

\begin{center}
  $X^{\omega-w} := \{ A\subseteq X : A\in \mathbf{WD}_{\omega}(X)\}$
\end{center}
topologized by the open sets $U^{ \omega-w}:=\{ A \in X ^{\omega-w} : A \cap U \neq \emptyset\}$ for each open subset $U$ of $X$. If we define $\eta^{\omega-w}_{X} :X \rightarrow X ^{\omega-w}$ by $\eta^{\omega-w}_{X} ( x ) = cl(\{ x \}) $, then we obtain a $\omega$-well-filterification \cite[Theorem~6.8]{XiaoquanXu}, which we call the standard $\omega$-well-filterification.
\begin{proposition}\label{continuous}
  Let $X,Y$ be two $\omega$-well-filtered space. If $f:X\rightarrow Y$ is a continuous map, then $f(\sup A)=\sup f(A)$ for any $A\in WD_{\omega}(X)$.
  \begin{proof}
    The topology of $X$ is denoted by $\tau$. For any $A\in WD_{\omega}(X)$, we know that $cl_{\tau}(A)=\da \sup A$ again by applying Corollary 6.11 in \cite{XiaoquanXu} since $X$ is $\omega$-well-filtered. Note that $f$ is monotone by virtue of being continuous. Then $f(\sup A)$ is an upper bound of $f(A)$. Let $t$ be any other upper bound of $f(A)$. This means that $A\subseteq f^{-1}(\da t)$. The fact that $f$ is continuous implies that $f^{-1}(\da t)$ is closed. We deduce that $\sup A\in cl_{\tau}(A)\subseteq f^{-1}(\da t)$. So $f(\sup A)\leq t$. Therefore, $f(\sup A)=\sup f(A)$.
\end{proof}
\end{proposition}
\begin{theorem}\label{well-filterification}
  Let $X$ be a core compactly generated space. Then the $\omega$-well-filterification of $X$ is again a core compactly generated space.
  \begin{proof}
  Let $(X^{\omega-w}, \eta_{X}^{\omega-w})$ be the standard $\omega$-well-filterification of $X$. We want to prove that $\eta_{X}^{\omega-w}:X\rightarrow (X^{\omega-w}, \mathcal{C}X^{\omega-w})$ is also continuous. To this end, let $U\in\mathcal{C}X^{\omega-w}$. Since $X$ is core compactly generated, it suffices to prove that $({\eta_{X}^{\omega-w}})^{-1}(U)\in \mathcal{C}X$.

   Assume $\rho:C \rightarrow X$ is a continuous function from a core compact space $C$. This means that $\eta_{X}^{\omega-w}\circ \rho$ is also continuous, which mainfests that $(\eta_{X}^{\omega-w}\circ \rho)^{-1}(U)$ is open in C. Now we can gain our desired result that $({\eta_{X}^{\omega-w}})^{-1}(U)\in \mathcal{C}X$. From Lemma \ref{well-filtered}, we know that $(X^{\omega-w},\mathcal{C}X^{\omega-w})$ is an $\omega$-well-filtered space. It follows that there is a unique continuous function $\overline{f}:X^{\omega-w}\rightarrow (X^{\omega-w}, \mathcal{C}X^{\omega-w})$ with $\overline{f}\circ\eta_{X}^{\omega-w}=\eta_{X}^{\omega-w}$.

    It is evident to see that $A=\sup_{a\in A}\da a$ for any $A\in \mathbf{WD}_{\omega}(X)$. From the definition of $\omega$-well-filtered determined subsets and the continuity of $\eta_{X}^{\omega-w}$, we have that $\{\da a\mid a\in A\}=\eta_{X}^{\omega-w}(A)$ is an $\omega$-well-filtered determined subset of $X^{\omega-w}$. In virtue of Proposition \ref{continuous}, we demonstrate that $\overline{f}(A)=\overline{f}(\sup _{a\in A}\da a)=\sup _{a\in A}\overline{f}(\da a)=\sup_{a\in A}\overline{f}\circ \eta_{X}^{\omega-w}(a)=\sup_{a\in A}\eta_{X}^{\omega-w}(a)=A$. It follows that $\overline{f}$ is an identity map. This suggests that the core compactly generated topology of $X^{\omega-w}$ is contained in $ \mathcal{O}(X^{\omega-w})$. Hence, $X^{\omega-w}$ is core compactly generated.
  \end{proof}
\end{theorem}
 From the definitions of the well-filterification and the $\omega$-well-filterification, one wonders whether the core compactly generated topology of a well-filtered space is again well-filtered and the well-filterification of a core compactly generated space is again core compactly generated. We end this section by leaving the above questions open.
\section{The core compactly generated topology of monotone convergence spaces}
In this section, we discuss the core compactly generated topology of monotone convergence spaces. We show that the core compactly generated topology of a monotone convergence space is again a monotone convergence space and the $D$-completion of a core compactly generated space is again core compactly generated.

Recall that a monotone convergence space $Y$ together with a topological embedding
$j:X \rightarrow Y$ with $j(X)$ a $d$-dense subset of $Y$ is called a $D$-completion of the $T_{0}$ space $X$. A subset $A$ of a $T_{0}$ space $X$ is called tapered if for any continuous function $f : X \rightarrow  Y$ mapping into a monotone convergence space $Y$, $ \sup f(A)$ always exists in $Y$.

A standard construction for the $D$-completion of a $T_{0}$ space $X$ is to set

\begin{center}
  $X^{d} := \{ A\subseteq X : A ~\mathrm{is} ~\mathrm{closed} ~\mathrm{and} ~\mathrm{tarped}\}$
\end{center}
topologized by open sets $U^{ d}:=\{ A \in X ^{d} : A \cap U \neq \emptyset\}$ for each open subset $U$ of $X$. If we define $\eta^{d}_{X} :X \rightarrow X ^{d}$ by $\eta^{d}_{X} ( x ) = cl(\{ x \}) $, then we obtain a $D$-completion \cite[Theorem 3.10]{zhongxizhang}, which we call the standard $D$-completion.

The following theorem is an immediate consequence of Lemma \ref{Scott}.
\begin{theorem}\label{monotone}
  Let $X$ be a monotone convergence space. Then the core compactly generated topology of $X$ is also monotone convergence.
\end{theorem}
With the help of Theorem \ref{monotone}, similar to the proof of Theorem \ref{well-filterification}, the next theorem follows directly.
\begin{theorem}
 Let $X$ be a core compactly generated space. Then the $D$-completion of $X$ is again a core compactly generated space.
\end{theorem}
The following lemma is a generation of \cite[Lemma 8.1]{Escardo}. We would mention the following:
\begin{lemma}\label{same}
  Let $f:X\rightarrow Y$ be continuous, where $X$ is a well-filtered space. If $Q$ is the intersection of a filtered family $\mathcal{K}$ of compact saturated subsets of $X$, then $\bigcap_{K\in \mathcal{K}}\ua f(K)=\ua f(Q)$.
  \begin{proof}
    It is easy to see that $\ua f(Q)\subseteq \bigcap_{K\in \mathcal{K}}\ua f(K)$. Note that $\ua f(Q)=\bigcap_{U\in \mathcal{U}}U$, where $\mathcal{U}=\{U\in \mathcal{O}(X)\mid f(Q)\subseteq U\}$. So it suffices to prove that $\bigcap_{K\in \mathcal{K}}\ua f(K)\subseteq U$ for any $U\in \mathcal{U}$. Let $U\in \mathcal{U}$. Then $Q=\bigcap_{K\in \mathcal{K}}K\subseteq f^{-1}(U)$. From the continuity of $f$, we can get that $f^{-1}(U)$ is open. The assumption that $X$ is well-filtered implies that there exists $K\in \mathcal{K}$ such that $K\subseteq f^{-1}(U)$. In other words, $\ua f(K)\subseteq U$.
  \end{proof}
\end{lemma}

From Lemma \ref{Scott}, we obtain the result that the core compactly generated topology of a monotone convergence space is contained in the Scott topology. Naturally, we wonder on which class of topological spaces, the core compactly generated topologies happen to be the Scott topologies determined by the specialization orders on those spaces. Now we show that the Smyth power spaces of well-filtered spaces satisfy the above property.

\begin{lemma}\label{coincide}
  Let $X$ be a well-filtered space. Then the topology determined by the continuous functions from locally compact sober spaces agrees with the core compactly generated topology.
  \begin{proof}
In the light of \cite[Theorem~3.1]{LWX}, we have that a core compact well-filtered space is locally compact sober. With the help of well-filterification, we can deduce the lemma by a similar proof of \cite[Lemma~8.2]{Escardo}.
  \end{proof}
\end{lemma}
\begin{theorem}
  Let $X$ be a well-filtered space. Then the core compactly generated topology of the Smyth power space $P_{S}(X)$ of $X$ equals the Scott topology under its specialization order.
  \begin{proof}
 By applying \cite[Theorem~5]{Xiaoquan22}, we know that $P_{S}(X)$ is well-filtered. Thus $P_{S}(X)$ is a monotone convergence space. It follows that the core compactly generated topology of $P_{S}(X)$ is contained in the Scott topology. Conversely, to this end, let $\mathcal{U}\in \sigma(P_{S}(X))$. Due to Lemma \ref{coincide}, we only need to prove that $\mathcal{U}$ is contained in the topology determined by the continuous functions from locally compact sober spaces.

 Assume $\rho:C\rightarrow X$ be a continuous function, where $C$ is a locally compact sober space. It remains to prove that $\rho^{-1}(\mathcal{U})$ is open. Now let $c\in \rho^{-1}(\mathcal{U})$. Since $C$ is locally compact, we have that $\ua c=\bigcap _{K\in \mathcal{R}}K$, where $\mathcal{R}=\{K\in Q(X)\mid c\in int(K)\}$. Owing to \cite[Lemma 8.1]{Escardo}, we can conclude that $\ua \rho(c)=\bigcap_{K\in \mathcal{R}}\ua \rho(K)\subseteq \mathcal{U}$.

 Note that $\rho(K)$ is compact in $P_{S}(X)$. Then $\bigcup \rho(K)$ is a compact subset of $X$.

 We claim that $\rho(c)=\bigcap_{K\in \mathcal{R}}(\bigcup\rho(K))$.

 For any $x\in \rho(c)$, it turns out that $\ua x\in \ua \rho(c)=\bigcap_{K\in \mathcal{R}}\ua \rho(K)$. This means that there is $k\in K$ such that $\ua x\subseteq \rho(k)$ for any $K\in \mathcal{R}$. Therefore, $x\in \rho(k)\subseteq \bigcup\rho(K)$ for any $K\in \mathcal{R}$, which is equivalent to saying that $x\in \bigcap_{ K\in \mathcal{R}}(\bigcup\rho(K))$.

 Conversely, let $x\in \bigcap_{K\in \mathcal{R}}(\bigcup\rho(K))$. This implies that $x\in \bigcup \rho(K)$ for any $K\in \mathcal{R}$, which guarantees the existence of $k\in K$ such that $x\in \rho(k)$ for any $K\in \mathcal{R}$. It follows that $\ua x\in \bigcap_{K\in \mathcal{R}}\ua \rho(K)=\ua \rho(c)$. So we have $x\in \rho(c)$.

 As we shall see that $\rho(c)=\bigcap_{K\in \mathcal{R}}(\bigcup\rho(K))\in \mathcal{U}$ and $\mathcal{R}$ is a directed subset of $P_{S}(X)$ in its specialization order, which yield that there exits $K\in \mathcal{R}$ such that $\bigcup\rho(K)\in \mathcal{U}$. The fact that $\mathcal{U}$ is an upper set suggests that $\rho(K)\subseteq \mathcal{U}$. We conclude that $c\in int(K)\subseteq K\subseteq \rho^{-1}(\mathcal{U})$. Now we can gain our desired result.
  \end{proof}
\end{theorem}
By applying the above lemma, we can gain the following corollaries which have been proved in \cite{Schalk} and \cite{XY}, respectively.
\begin{theorem}(\cite{Schalk})
If $X$ is a locally compact sober space, then the upper Vietoris topology and the Scott topology on $Q(X)$ coincide.
\end{theorem}
\begin{theorem}(\cite{XY})
If $X$ is a well-filtered space and its Smyth power space $P_{S}(X)$ is first-countable, then the upper Vietoris topology
agrees with the Scott topology on $Q(X)$.
\end{theorem}
\vskip 0.2cm
\section{The product of two core compactly generated spaces}

If $X,Y$ are two core compactly generated spaces, then the topological product $X\times Y$ need not be a core compactly generated topology, in general. One always has $\mathcal{O}(X\times Y)\subseteq \mathcal{C}(X\times Y)$, but the containment may be proper \cite[Example~5.1]{Escardo}. In this section, we obtain a necessary and sufficient condition for $\mathcal{O}(X\times Y)= \mathcal{C}(X\times Y)$.

\begin{theorem}
  Let $X$ be a core compactly generated space. Then the following statements are equivalent.

  (1) $X$ is core compact.

  (2) For every core compactly generated space $Y$, one has $\mathcal{O}(X\times Y)=\mathcal{C}(X\times Y)$.
  \begin{proof}
    $(1)\Rightarrow(2)$ follows from \cite[Theorem~5.4]{Escardo}.

    $(2)\Rightarrow (1)$: Let $Y=\Sigma (\mathcal{O}(X))$. Then $Y$ is a core compactly generated space by  \cite[Lemma~4.6]{Escardo}. It follows that $\mathcal{O}(X\times Y)=\mathcal{C}(X\times Y)$ since $(2)$ holds. Due to \cite[Exersice~5.2.7]{GL}, we can know that $X$ is core compact iff $G=\{(x,U)\in X\times \mathcal{O}(X)\mid x\in U\}$ is open in $X\times \Sigma \mathcal{O}(X)$. So it suffices to prove that $G$ is core compactly generated open.

    Assume $\rho: C\rightarrow X$ be continuous, where $C$ is a core compact space. It remains to prove that $\rho^{-1}(G)$ is open in $C$. For any $c\in \rho^{-1}(G)$, we have that $\rho(c)\in G$. We write $\rho(c)=(x_{c},U_{c})$. Then $x_{c}\in U_{c}$. Note that $x_{c}=P_{X}\circ \rho(c)$, where $P_{X}: X\times \Sigma \mathcal{O}(X)\rightarrow X$ is the projection map. This manifests that $c\in (P_{X}\circ \rho)^{-1}(U_{c})$, $(P_{X}\circ \rho)^{-1}(U_{c})$ is open.  In the light of the core compactness of $C$, we construct inductively a decreasing family of open subsets with
    \begin{center}
      $x\in V\ll \cdots \ll V_{n+1}\ll V_{n}\ll\cdots\ll V_{1}= (P_{X}\circ \rho)^{-1}(U_{c})$.
    \end{center}

    We define $A_{n}=\{U\in \mathcal{O}(X)\mid P_{X}\circ \rho(V_{n})\subseteq U\}$ for any $n\in \mathbb{N}$ and $R=\bigcup_{n\in \mathbb{N}}A_{n}$.

    We claim that $R$ is Scott open.

    Clearly, $R$ is an upper set. Let $(U_{i})_{i\in I}$ be a directed subset of $\mathcal{O}(X)$ with $\bigcup_{i\in I}U_{i}\in R$. Then there is $n\in \mathbb{N}$ such that $\bigcup_{i\in I}U_{i}\in A_{n}$, in other words, $P_{X}\circ \rho(V_{n})\subseteq \bigcup_{i\in I}U_{i}$. The fact that $V_{n+1}\ll V_{n}$ suggests that there is $i\in I$ such that $P_{X}\circ \rho(V_{n+1})\subseteq U_{i}$, that is, $U_{i}\in A_{n+1}\subseteq R$.

    Additionally $P_{\mathcal{O}(X)}\circ \rho(c)=U_{c}\in R$, where $P_{\mathcal{O}(X)}:X\times \mathcal{O}(X)\rightarrow \mathcal{O}(X)$ is the projection map. This yields that there exists an open set $U$ of $X$ such that $c\in W\subseteq (P_{\mathcal{O}(X)}\circ \rho)^{-1}(R)$.

    It suffices to show that $c\in W\cap V\subseteq \rho^{-1}(G)$.

    To this end, let $a\in W\cap V$. Then $P_{\mathcal{O}(X)}\circ \rho(a)\in R$. The construction of $R$ guarantees the existence of $n\in \mathbb{N}$ with $P_{X}\circ \rho(V_{n})\subseteq P_{\mathcal{O}(X)}\circ \rho(a)$. As we shall see $a\in V\subseteq V_{n}$, which means that $P_{X}\circ \rho(a)\subseteq P_{X}\circ \rho(V_{n})\subseteq P_{\mathcal{O}(X)}\circ \rho(a)$. Hence, $\rho(a)\in G$.
\end{proof}
\end{theorem}
\section{Reference}
\bibliographystyle{plain}

\begin{thebibliography}{10}

\bibitem{Battenbeld} O. Battenbeld, M. Schr$\mathrm{\ddot{o}}$der, A. Simpson, Compactly generated domain theory, Mathematical Structures in Computer Science. 16 (2006) 141-161.
  \bibitem{Day} B. Day, A reflection theorem for closed categories, Journal of Pure and Applied Algebra. 2 (1972) 1-11.
\bibitem{RE} R. Engelking, General Topology, Polish Scientific Publishers, 1989.
\bibitem{Escardo}M. Escard$\mathrm{\acute{o}}$, J. Lawson, A. Simpson, Comparing Cartesian closed categories of (core) compactly generated spaces, Toplogy and its applications. 143 (2004) 105-145.
\bibitem{Ershov}Y. Ershov, The $d$-rank of a topological space, Algebra and Logic. 56 (2017) 98-107.
\bibitem{M. Escardo} M. Escardo, Compactly generated Hausdorff locales, Annals of Pure and Applied Logic. 137 (2006) 147-163.

\bibitem{GHKLMS} G. Gierz, K. Hofmann, K. Keimel, J. Lawson, M. Mislove, D. Scott, Continuous Lattices and Domains, Cambridge University Press, 2003.

\bibitem{GL} J. Goubault-Larrecq, Non-Hausdorff Topology and Domain Theory, New Mathematical Monographs, Cambridge University Press, 2013.
\bibitem{Gruenhage} G. Gruenhage, T. Streicher, Quotients of countably based spaces are not closed under sobrification, Mathematical Structures in Computer Science. 16 (2006) 223-229.
\bibitem{Jia} X. Jia, Meet-continuity and locally compact sober dcpos, PhD thesis, University of Birmingham, 2018.
\bibitem{Keimel} K. Keimel, J. Lawson, $D$-completions and the $d$-topology, Annals of Pure and Applied Logic. 159 (2007) 292-306.
\bibitem{LWX} J. Lawson, G. Wu, X. Xi, Well-filtered spaces, compactness, and the lower topology, Houston Journal of Mathematics. 46 (2020) 283-294.
\bibitem{Markowsky}G. Markowsky, Chain-complete posets and directed sets with applications, Algebra Universalis. 6 (1976) 53-68.

\bibitem{Miao1} H. Miao, Z. Yuan, Q. Li, A discussion of well-filteredness and sobriety, Topology and its Applications. 291 (2020) 107450.

\bibitem{Miao2} H. Miao, Q. Li, D. Zhao, On two problems about sobriety of topological spaces, Topology and its Applications. 295 (2021) 107667.

\bibitem{Miao3} H. Miao, L. Wang, Q. Li, D-completion, well-filterification and sobrification, arXiv: 2101.04894.


\bibitem{SXXZ} C. Shen, X. Xi, X. Xu, D. Zhao, On well-filtered reflections of $T_{0}$ space, Topology and its Applications. 267 (2019) 106869.
\bibitem{Schalk} A. Schalk, Algebras for Generalized Power Constructions, PhD Thesis, Technische Hochschule Darmstadt, 1993.

\bibitem{Vogt} R. Vogt, Convenient categories of topological spaces for homotopy theory, Archiv der Mathematik (Basel). 22 (1971) 545-555.
\bibitem{Wyler} O. Wyler, Dedekind complete posets and Scott topologies, in: B. Banaschewsk, R.-E. Hoffmann (Eds.), Continuous Lattices Proceedings, Bremen 1979, vol. 871, Springer, Berlin, Heidelberg, 1981, pp. 384-389.
\bibitem{WXXZ} G. Wu, X. Xi, X. Xu, D. Zhao, Existence of well-filterifications of $T_{0}$ topological spaces, Topology and its Applications. 270 (2019) 107044.
\bibitem{XSXZ} X. Xu, C. Shen, X. Xi, D. Zhao, On $T_{0}$ spaces determined by well-filtered spaces, Topology and its Applications. 282 (2020) 107323.
\bibitem{XiaoquanXu} X. Xu, C. Shen, X. Xi, D. Zhao, First countability, $\omega$-well-filtered spaces and reflections, Topology and its Applications. 279 (2020) 107255.
\bibitem{xilawson} X. Xi, J. Lawson, On well-filtered spaces and ordered sets, Topology and its Applications. 228 (2017) 139-144.
\bibitem{Xiaoquan22} X. Xu, X. Xi, D. Zhao, A complete Heyting algebra whose Scott space is non-sober, Fundamenta Mathematicae, 252 (2021) 315-323.
\bibitem{XY} X. Xu, Z. Yang, Coincidence of the upper Vietoris topology and the Scott topology, Topology and its Applications, 288 (2021) 107480.
\bibitem{zhongxizhang}Z. Zhang, Q. Li, A direct characterization of monotone convergence space completion, Topology and its Applications. 230 (2017) 99-104.
\bibitem{ZF}D. Zhao, T. Fan, Dcpo-completion of posets, Theoretical Computer Science. 411 (2010) 2167-2173.
\end{thebibliography}

\end{document}